\documentclass[12pt, twoside, leqno]{article}
\usepackage{amsmath,amsthm}
\usepackage{amssymb,latexsym}
\usepackage{enumerate}
\usepackage[T1]{fontenc}
\usepackage{array}
\usepackage[all]{xypic}
\usepackage{fancyhdr}

\newtheorem{thm}{Theorem}[section]
\newtheorem{cor}[thm]{Corollary}
\newtheorem{lem}[thm]{Lemma}
\newtheorem{prop}[thm]{Proposition}

\theoremstyle{definition}

\numberwithin{equation}{section}

 1

\def\CCC{{\mathbb{C}}}

\def\NN{{\mathbb{N}}}
\def\MM{{\mathbb{M}}}

\def\ZZ{{\mathbb{Z}}}

\def\RR{{\mathbb{R}}}

\def\CH{{\cal H}}

\def\CN{{\cal N}}

\def\CL{{\cal L}}

\def\CP{{\cal P}}

\def\CZ{{\cal Z}}
\def\CM{{\cal M}}

\def\CZ{{\cal Z}}

\def\epv {{$\mbox{}$\hfill ${\Box}$\vspace*{1.5ex} }}

\newcommand{\ra}{\rightarrow}

\renewcommand*{\dim}{\mathrm{dim}}

\def\disc{\textnormal{disc}}

\def\ov{\overline}
\def\wt{\widetilde}

\def\epv {{$\mbox{}$\hfill ${\Box}$\vspace*{1.5ex} }}

\def\ov#1{\overline{#1}}

\begin{document}

%%%%% To ease editing, for IMPAN journals add:

\baselineskip=17pt

%%%%%%%%%%%%%%%%

\title{Generalized Shemesh criterion, common invariant subspaces and irreducible completely positive superoperators}

\author{Andrzej Jamio{\l}kowski (Toru\'n) and Grzegorz Pastuszak (Warszawa)
}

\date{}

\maketitle

%% Classification and key words; note that the 2010 classification is used:

\renewcommand{\thefootnote}{}



\renewcommand{\thefootnote}{\arabic{footnote}}
\setcounter{footnote}{0}

%%%%%%%%

\begin{abstract}
Assume that $A_{1},...,A_{s}$ are complex $n\times n$ matrices. We give a computable criterion for existence of a common eigenvector of $A_{i}$ which generalize the result of D. Shemesh established for two matrices. We use this criterion to prove some necessary and sufficient condition for $A_{i}$ to have a common invariant subspace of dimension $d$, $2\leq d<n$, if every $A_{i}$ has pairwise different eigenvalues. Finally, we observe that the set of all matrices having multiple eigevalues has Lebesgue measure $0$ and thus the condition is sufficient in practical applications. Being motivated by quantum information theory, we give a flavour of such applications for irreducible completely positive superoperators.
\end{abstract}

\section{Introduction and the main results}

Quantum theory, in its nonrelativistic formulation, is built on the theory of Hilbert spaces and operators. Assume that $\CH$ is a fixed Hilbert space associated with a given quantum system $S$. We denote by $B(\CH)$ the set of all linear continuous operators on $\CH$. Then the set of states of the system $S$ is, by definition, represented by all semipositive elements of $B(\CH)$ with trace equal to one. This set of states will be denoted by $S(\CH)$. 

In the beginning of seventies it appeared that some natural questions connected with fundamentals of quantum mechanics (more precisely, with the theory of open quantum systems) lead to investigations of linear maps in a real Banach space of self-adjoint operators on a fixed Hilbert space (cf. e.g. \cite{HZ}, \cite{Ja}). The concept of a Banach space with the partial order defined by a specific cone, namely, the cone of positive semidefinite operators, constitutes a basic idea in the description of open quantum systems and in the quantum information theory (\cite{HZ}, cf. also the last section of this paper). 

In this paper we consider finite dimensional Hilbert spaces - such spaces are the current focus in quantum computing and quantum information theory for experimental reasons. This means we assume $\CH\cong\CCC^{n}$ and $B(\CH)\cong\MM_{n}(\CCC)$. It should be stressed that to describe all possible changes of quantum states one has to consider some specific linear operators in $B(\CH)$. Very often, at least in physical literature, they are called \textit{superoperators}. The general form of such maps is well known. Namely, for a given superoperator $\Phi:B(\CH)\ra B(\CH)$, $\dim\CH<\infty$, there always exists an operator-sum presentation given by $$(*)\quad\Phi(X)=\sum_{i=1}^{s}A_{i}XB_{i},$$ where $A_{i},B_{i}$ are elements of $B(\CH)$, for all $i=1,...,s$. A particular class of such maps, the so-called \textit{completely positive maps}  (or in physical terminology \textit{quantum operations} or \textit{quantum channels}) play a prominent role in formulations of evolution of open quantum systems and in the theory of quantum measurements.

It is obvious from the above considerations that properties of superoperators $\Phi:B(\CH)\ra B(\CH)$ are connected with properties of the sets $\{A_{1},...,A_{s}\}$ and $\{B_{1},...,B_{s}\}$ of matrices. In particular, in case of completely positive maps, we have $B_{i}=A_{i}^{*}$ for $i=1,...,s$, where $A_{i}^{*}$ denotes matrix adjoint to $A_{i}$, and thus quantum channels are described by maps of the form $$(**)\quad\Phi(X)=\sum_{i=1}^{s}A_{i}XA_{i}^{*}.$$

The matrices $A_{i}$ occuring in the expression above are called the \textit{Kraus coefficients} of $\Phi$.

Motivated by the main results of \cite{Fa}, we discuss in the paper some properties of a completely positive map $\Phi$ of the form $(**)$ in terms of common invariant subspaces of its Kraus coefficients $A_{i}$. Namely, it turns out that $\Phi$ is \textit{irreducible} (see \cite{Fa} for the definitions) if and only if the matrices $A_{i}$ do not have a nontrivial common invariant subspace. 

This paper is devoted to give a computable criterion for a completely positive superoperator $\Phi$ to be irreducible. To establish the criterion we give at first a computable condition for the matrices $A_{1},...,A_{s}$ to have a common invariant subspace of a fixed dimension $d$.

In what follows, we consequently assume that $B(\CH)\cong\MM_{n}(\CCC)$ is the vector space of all $n\times n$ complex matrices.

Assume that $x\in\CCC^{n}$, $x\neq 0$. We say that $x$ is a \textit{common eigenvector} of $A_{1},...,A_{s}\in\MM_{n}(\CCC)$ if and only if $x$ is an eigenvector of every $A_{i}$, that is, $A_{i}x=\mu_{i} x$ for some $\mu_{i}\in\CCC$. 

Assume that $W$ is a subspace of $\CCC^{n}$. We say that $W$ is a \textit{common invariant subspace} of $A_{1},...,A_{s}\in\MM_{n}(\CCC)$ if and only if $W$ is $A_{i}$-invariant for all $i=1,...,s$, that is, $A_{i}w\in W$ for all $w\in W$. It is clear that $A_{i}$ have a common eigenvector if and only if $A_{i}$ have a common invariant subspace of dimension $1$.

Assume that $A,B\in\MM_{n}(\CCC)$. We denote by $[A,B]:=AB-BA$ the \textit{commutator} of $A$ and $B$, and by $\ker A:=\{v\in\CCC^{n}|Av=0\}$ the \textit{kernel} of $A$. 

In \cite{Sh} D. Shemesh proved the following criterion for existence of common eigenvector of two complex matrices $A$ and $B$. 

\begin{thm} \textnormal{(\cite[Theorem 3.1]{Sh})} Assume that $A,B\in\MM_{n}(\CCC)$ and $$\CN(A,B):=\bigcap_{k,l=1}^{n-1}\ker[A^{k},B^{l}].$$ Then $A$ and $B$ have a common eigenvector if and only if $\CN(A,B)\neq 0$.
\end{thm}

It it observed in \cite{Sh} that $\CN(A,B)=\ker K$, where 

$$K=\sum_{k,l=1}^{n-1}[A^{k},B^{l}]^{*}[A^{k},B^{l}]$$

and $X^{*}$ denotes matrix adjoint to $X$. It follows that the condition $$\CN(A,B)=\bigcap_{k,l=1}^{n-1}\ker[A^{k},B^{l}]\neq 0$$ is \textit{computable}, i.e. it can be verified by finite and deterministic algorithm.

In view of Theorem 1.1 it is natural to consider the following problem: \textit{Is there a computable condition verifing the existence of common invariant subspace of dimension $d$ of $s$ complex $n\times n$ matrices $A_{1},...,A_{s}$, where $1\leq d<n$ and $s\geq 2$}?

Partial solutions to this problem are given in \cite{AGI}, \cite{AI}, \cite{GI} and \cite{Ts} where it is mainly assumed that $s=2$. In \cite{AGI} and \cite{AI} the authors study the case when algebra generated by two complex matrices is semisimple and use the concept of a \textit{standard polynomial}, see \cite{Pi}, Section 20.4. In \cite{GI} and \cite{Ts} the authors reduce the general problem to the question of existence of common eigenvector of suitable \textit{compound matrices}, see \cite{La}, Chapter I.

The general version of the problem with arbitrary $d$ and $s$ is solved completely in \cite{AP} where some techniques of Gr{\"o}bner bases theory and algebraic geometry are used.

%This problem is studied in \cite{AGI}, \cite{AI}, \cite{GI} and \cite{Ts}, and each of these papers gives to it a partial solution. In \cite{Ts} the problem was completely solved theoretically, but the proposed algorithm involved steps that were hard to verify in a practical manner. 

%A major breakthrough has been made in \cite{AP} where the authors gave a computable condition for $A_{1},...,A_{s}\in\MM_{n}(\CCC)$ to have a common invariant subspace of dimension $d$, where $1\leq d<n$ and $s\geq 2$. In \cite{AP} completely new approach to the problem of common invariant subspaces is taken. Namely, the authors use some methods of Grobner bases theory and algebraic geometry in their solution, whereas in earlier papers mainly basic linear algebra is used. 

In this paper we present a different approach to the problem discussed. Namely, for $d=1$ we prove the following computable generalization of the Shemesh criterion.

\begin{thm} Assume that $A_{1},...,A_{s}\in\MM_{n}(\CCC)$ and $$\CM(A_{1},...,A_{s}):=\bigcap_{\begin{smallmatrix} k_{i},l_{j}\geq 0 \\ k_{1}+k_{2}+...+k_{s}\neq 0\\l_{1}+l_{2}+...+l_{s}\neq 0\end{smallmatrix}}^{n-1}\ker [A_{1}^{k_{1}}...A_{s}^{k_{s}},A_{1}^{l_{1}}...A_{s}^{l_{s}}].$$ 

{\rm (1)} Matrices $A_{i}$ have a common eigenvector if and only if $$\CM(A_{1},...,A_{s})\neq 0.$$

{\rm (2)} We have $\CM(A_{1},...,A_{s})= \ker K$ where
 
$$K=\sum_{\begin{smallmatrix} k_{i},l_{j}\geq 0 \\ k_{1}+k_{2}+...+k_{s}\neq 0\\l_{1}+l_{2}+...+l_{s}\neq 0\end{smallmatrix}}^{n-1}[A_{1}^{k_{1}}...A_{s}^{k_{s}},A_{1}^{l_{1}}...A_{s}^{l_{s}}]^{*}[A_{1}^{k_{1}}...A_{s}^{k_{s}},A_{1}^{l_{1}}...A_{s}^{l_{s}}].$$
\end{thm}

If $d>1$ and each of $A_{1},...,A_{s}\in\MM_{n}(\CCC)$ has pairwise different eigenvalues, we apply Theorem 1.2 to establish a computable criterion for existence of common invariant subspace of $A_{i}$ of dimension $d$. In this approach we make use of some methods presented in \cite{GI} and \cite{Ts}. 

Furthermore, we observe that the set of all matrices having at least one multiple eigevalue is Lebesgue-measureble, and its measure is equal to zero. This yields that when matrices $A_{1},...,A_{s}\in\MM_{n}(\CCC)$ are random it is reasonble to expect that non of them has multiple eigenvalues. Consequently, our criterion may be considered as sufficient for practical applications. 

In the final section of the paper we give a flavor of such applications presenting a computable method of verifing whether a completely positive superoperator defined by its Kraus coefficients is irreducible. 

\section{The generalized Shemesh criterion}

This section is devoted to the detailed proof of Theorem 1.2 which we shall call \textit{the generalized Shemesh criterion}.

We start with a simple fact from linear algebra.

\begin{prop} Assume that $A\in\MM_{n}(\CCC)$ and $X\neq 0$ is $A$-invariant. Then there is an eigenvector $w\in\CCC^{n}$ of $A$ such that $w\in X$.
\end{prop}

{\bf Proof.} Assume that $\dim X=s\geq 1$ and $X=\langle x_{1},...,x_{s}\rangle$. The Steinitz exchange lemma implies that there are vectors $y_{s+1},...,y_{n}\in\CCC^{n}$ such that $\CCC^{n}=\langle x_{1},...,x_{s},y_{s+1},...,y_{n}\rangle$. Assume that $B$ is an invertible matrix of the form $$B=\left[\begin{matrix}x_{1}&x_{2}&\ldots & x_{s} &y_{s+1}&\ldots &y_{n}\end{matrix}\right].$$ Since $X$ is $A$-invariant, there is a matrix $$\wt{A}=\left[\begin{matrix}A_{1}& A_{2}\\0&A_{3}\end{matrix}\right],$$ where $A_{1}\in\MM_{s}(\CCC)$, such that $A=B\wt{A}B^{-1}$. Clearly there is $v=(v_{1},...,v_{s})^{tr}\in\CCC^{s}$ with $A_{1}v=\mu v$ for some $\mu\in\CCC$. This implies that $\wt{A}\wt{v}=\mu\wt{v}$ for $\wt{v}=(v_{1},...,v_{s},\underbrace{0,...,0}_{n-s})^{tr}$ and hence $A(B\wt{v})=\mu(B\wt{v})$ and $B\wt{v}\in X$. Thus $w:=B\wt{v}$ satisfies the required conditions.\epv

Assume that $A\in\MM_{n}(\CCC)$ and $\mu\in\CCC$ is an eigenvalue of $A$. We denote by $V_{\mu}^{A}$ the eigenspace associated to $\mu$, that is, $$V_{\mu}^{A}:=\{v\in\CCC^{n}|Av=\mu v\}.$$

The following lemma is commonly known (see \cite{Dr}, \cite{DDG}, \cite{McC} for more general versions). It states that pairwise commuting square complex matrices $A_{1},...,A_{s}\in\MM_{n}(\CCC)$ have a common eigenvector. 

\begin{lem} Assume that $A_{1},...,A_{s}\in\MM_{n}(\CCC)$ and $A_{i}A_{j}=A_{j}A_{i}$ for all $i,j=1,...,s$. Then there is a common eigenvector of $A_{i}$.
\end{lem}

{\bf Proof.} Assume that $s=2$ and $A=A_{1}$, $B=A_{2}$. Since $A$ is a complex matrix, $A$ has at least one eigenvalue, say $\mu\in\CCC$. If $x\in V_{\mu}^{A}$, then $$A(Bx)=(AB)x=(BA)x=B(Ax)=B(\mu x)=\mu(Bx)$$ hence $V_{\mu}^{A}$ is $B$-invariant. It follows from Proposition 2.1 that there is an eigenvector $y$ of $B$ such that $y\in V_{\mu}^{A}$. Thus $y$ is a common eigenvector of $A$ and $B$.

Assume now that $s\geq 3$ and $x\in\CCC^{n}$ is a common eigenvector of $A_{1},...,A_{s-1}$ such that $A_{1}x=\mu_{1}x$,..., $A_{s-1}x=\mu_{s-1}x$, for some $\mu_{1},...,\mu_{s-1}\in\CCC$, and hence $x\in V_{\mu_{1}}^{A}\cap ...\cap V_{\mu_{s-1}}^{A}$. We will show that $A_{1},...,A_{s-1},A_{s}$ also have a common eigenvector.

Indeed, observe that $$A_{i}(A_{s}x)=(A_{i}A_{s})x=(A_{s}A_{i})x=A_{s}(A_{i}x)=A_{s}(\mu_{i}x)=\mu_{i}(A_{s}x)$$ and hence $V_{\mu_{i}}^{A}$ is $A_{s}$-invariant for any $i=1,...,s-1$. It follows that $V_{\mu_{1}}^{A}\cap ...\cap V_{\mu_{s-1}}^{A}$ is a nonzero $A_{s}$-invariant subspace and hence there is an eigenvector $y$ of $A_{s}$ such that $y\in V_{\mu_{1}}^{A}\cap ...\cap V_{\mu_{s-1}}^{A}$, by Proposition 2.1. Thus $y$ is is a common eigenvector of $A_{1},...,A_{s-1},A_{s}$, and the lemma follows by induction.\epv

\begin{cor} Assume that $A_{1},...,A_{s}\in\MM_{n}(\CCC)$, $X\neq 0$ is $A_{i}$-invariant and $(A_{i}A_{j})x=(A_{j}A_{i})x$ for any $x\in X$. Then there is a common eigenvector $u$ of $A_{i}$ such that $u\in X$.
\end{cor}

{\bf Proof.} Assume that $\dim X=t\geq 1$ and $X=\langle x_{1},...,x_{t}\rangle$. Steinitz exchange lemma implies that $\CCC^{n}=\langle x_{1},...,x_{t},y_{t+1},...,y_{n}\rangle$ for some $y_{t+1},...,y_{n}\in\CCC^{n}$. We set $$B=\left[\begin{matrix}x_{1}&x_{2}&\ldots & x_{t} &y_{t+1}&\ldots &y_{n}\end{matrix}\right].$$ Since $X$ is $A_{i}$-invariant, there are matrices $$\wt{A_{i}}=\left[\begin{matrix}A^{i}_{1}& A^{i}_{2}\\0&A^{i}_{3}\end{matrix}\right],$$ where $A^{i}_{1}\in\MM_{t}(\CCC)$, such that $A_{i}=B\wt{A_{i}}B^{-1}$. We will show that $A^{i}_{1}A^{j}_{1}=A^{j}_{1}A^{i}_{1}$ for any $i,j=1,...,s$.

Assume that $v=(v_{1},...,v_{t})^{tr}\in\CCC^{t}$ and $\wt{v}=(v_{1},...,v_{t},\underbrace{0,...,0}_{n-t})^{tr}\in\CCC^{n}$. Then clearly $B\wt{v}\in X$ and hence $$(A^{i}_{1}A^{j}_{1})v=(\wt{A_{i}}\wt{A_{j}})\wt{v}=(B^{-1}B)(\wt{A_{i}}\wt{A_{j}})(B^{-1}B)\wt{v}=B^{-1}(B(\wt{A_{i}}\wt{A_{j}})B^{-1})B\wt{v}=$$$$=B^{-1}(A_{i}A_{j})B\wt{v}=B^{-1}(A_{j}A_{i})B\wt{v}=B^{-1}(B(\wt{A_{j}}\wt{A_{i}})B^{-1})B\wt{v}=$$$$=(B^{-1}B)(\wt{A_{j}}\wt{A_{i}})(B^{-1}B)\wt{v}=(\wt{A_{j}}\wt{A_{i}})\wt{v}=(A^{j}_{1}A^{i}_{1})v$$ since $A_{i},A_{j}$ commute on $X$. This implies that $(A^{i}_{1}A^{j}_{1})e_{r}=(A^{j}_{1}A^{i}_{1})e_{r}$ for any vector $e_{r}$ of the standard basis of $\CCC^{t}$. Hence $A^{i}_{1}A^{j}_{1}=A^{j}_{1}A^{i}_{1}$ and it follows by Lemma 2.2 that $A_{1}^{i}$ have a common eigenvector $w=(w_{1},...,w_{t})^{tr}$. This implies that $B\wt{w}\in X$ is a common eigenvector of $A_{i}$, where $\wt{w}=(w_{1},...,w_{t},\underbrace{0,...,0}_{n-t})^{tr}$. Thus $u:=B\wt{w}$ satisfies the required conditions. \epv

Now we are ready to prove a generalization of \cite[Theorem 3.1]{Sh}.

\begin{thm} Assume that $A_{1},...,A_{s}\in\MM_{n}(\CCC)$ and $$\CM(A_{1},...,A_{s}):=\bigcap_{\begin{smallmatrix} k_{i},l_{j}\geq 0 \\ k_{1}+k_{2}+...+k_{s}\neq 0\\l_{1}+l_{2}+...+l_{s}\neq 0\end{smallmatrix}}^{\infty}\ker [A_{1}^{k_{1}}...A_{s}^{k_{s}},A_{1}^{l_{1}}...A_{s}^{l_{s}}].$$

{\rm (1)} The subspace $\CM(A_{1},...,A_{s})$ is $A_{i}$-invariant for any $i=1,...,s$.

{\rm (2)} Matrices $A_{i}$ have a common eigenvector if and only if $$\CM(A_{1},...,A_{s})\neq 0.$$
\end{thm}

{\bf Proof.} (1) If $\CM(A_{1},...,A_{s})=0$ then clearly $\CM(A_{1},...,A_{s})$ is $A_{i}$-invariant. Hence assume that $\CM(A_{1},...,A_{s})\neq 0$ and let $v\in\CM(A_{1},...,A_{s})$. We set $$\CP_{v}:=\{p(A_{1},...,A_{s})v|p\in\CCC\langle x_{1},...,x_{s}\rangle\},$$ where $\CCC\langle x_{1},...,x_{s}\rangle$ is a complex polynomial ring of $s$ non-commutative variables $x_{1},...,x_{s}$. We will show that $\CP_{v}\subseteq\CM(A_{1},...,A_{s})$ for any $v\in\CM(A_{1},...,A_{s})$. 

Assume that $k_{1},...,k_{s},l_{1},...,l_{s}\geq 0$. If $k_{2}+...+k_{s}\neq 0$ and $l_{1}+...+l_{s}\neq 0$, then
$$(A_{1}^{k_{1}}...A_{s}^{k_{s}})(A_{1}^{l_{1}}...A_{s}^{l_{s}})v=A_{1}^{k_{1}}[(A_{2}^{k_{2}}...A_{s}^{k_{s}})(A_{1}^{l_{1}}...A_{s}^{l_{s}})]v=$$$$=A_{1}^{k_{1}}[(A_{1}^{l_{1}}...A_{s}^{l_{s}})(A_{2}^{k_{2}}...A_{s}^{k_{s}})]v=A_{1}^{k_{1}+l_{1}}[(A_{2}^{l_{2}}...A_{s}^{l_{s}})(A_{2}^{k_{2}}...A_{s}^{k_{s}})]v$$ 

since $v\in\ker[A_{2}^{k_{2}}...A_{s}^{k_{s}},A_{1}^{l_{1}}...A_{s}^{l_{s}}]$.

If $k_{2}+...+k_{s}=0$ or $l_{1}+...+l_{s}=0$ then also
$$(A_{1}^{k_{1}}...A_{s}^{k_{s}})(A_{1}^{l_{1}}...A_{s}^{l_{s}})v=A_{1}^{k_{1}+l_{1}}[(A_{2}^{l_{2}}...A_{s}^{l_{s}})(A_{2}^{k_{2}}...A_{s}^{k_{s}})]v.$$

Therefore a simple induction argument implies that 
$$(A_{1}^{k_{1}}...A_{s}^{k_{s}})(A_{1}^{l_{1}}...A_{s}^{l_{s}})v=(A_{1}^{k_{1}+l_{1}}...A_{1}^{k_{s}+l_{s}})v$$ 

for any $k_{1},...,k_{s},l_{1},...,l_{s}\geq 0$. It follows that if $q\in\CCC\langle x_{1},...,x_{s}\rangle$, then $$q(A_{1}...A_{s})v=(\sum_{0\leq t_{i}\leq\deg(q)}a_{t_{1}...t_{s}}^{q}A_{1}^{t_{1}}...A_{s}^{t_{s}})v$$ for some $a_{t_{1}...t_{s}}^{q}\in\CCC$.

Moreover,
$$(A_{1}^{k_{1}}...A_{s}^{k_{s}})(A_{1}^{l_{1}}...A_{s}^{l_{s}})q(A_{1}...A_{s})v=$$$$=(A_{1}^{k_{1}}...A_{s}^{k_{s}})(A_{1}^{l_{1}}...A_{s}^{l_{s}})(\sum_{0\leq t_{i}\leq\deg(q)}a_{t_{1}...t_{s}}^{q}A_{1}^{t_{1}}...A_{s}^{t_{s}})v=$$$$=\sum_{0\leq t_{i}\leq\deg(q)}a_{t_{1}...t_{s}}^{q}(A_{1}^{k_{1}}...A_{s}^{k_{s}})(A_{1}^{l_{1}}...A_{s}^{l_{s}})A_{1}^{t_{1}}...A_{s}^{t_{s}}v=$$$$=\sum_{0\leq t_{i}\leq\deg(q)}a_{t_{1}...t_{s}}^{q}A_{1}^{t_{1}+k_{1}+l_{1}}...A_{s}^{t_{s}+k_{s}+l_{s}}v=$$$$=\sum_{0\leq t_{i}\leq\deg(q)}a_{t_{1}...t_{s}}^{q}(A_{1}^{l_{1}}...A_{s}^{l_{s}})(A_{1}^{k_{1}}...A_{s}^{k_{s}})A_{1}^{t_{1}}...A_{s}^{t_{s}}v=$$$$=(A_{1}^{l_{1}}...A_{s}^{l_{s}})(A_{1}^{k_{1}}...A_{s}^{k_{s}})(\sum_{0\leq t_{i}\leq\deg(q)}a_{t_{1}...t_{s}}^{q}A_{1}^{t_{1}}...A_{s}^{t_{s}})v=$$$$=(A_{1}^{l_{1}}...A_{s}^{l_{s}})(A_{1}^{k_{1}}...A_{s}^{k_{s}})q(A_{1}...A_{s})v$$ and thus $q(A_{1},...,A_{s})v\in\ker[A_{1}^{k_{1}}...A_{s}^{k_{s}},A_{1}^{l_{1}}...A_{s}^{l_{s}}]$ for any $q\in\CCC\langle x_{1},...,x_{s}\rangle$ and $k_{1},...,k_{s},l_{1},...,l_{s}\geq 0$. 

This implies that $\CP_{v}\subseteq\CM(A_{1},...,A_{s})$ for any $v\in\CM(A_{1},...,A_{s})$ and hence 
$$\CM(A_{1},...,A_{s})=\bigcup_{v\in\CM(A_{1},...,A_{s})}\CP_{v}.$$ 

This shows that the subspace $\CM(A_{1},...,A_{s})$ is $A_{i}$-invariant since clearly $\CP_{v}$ is $A_{i}$-invariant, for any $i=1,...,s$. 

(2) $\Rightarrow$ Assume that $x\in\CCC^{n}$ is a common eigenvector of $A_{i}$ such that $A_{i}x=\mu_{i}x$ for some $\mu_{i}\in\CCC$. Then 
 $$(A_{1}^{k_{1}}...A_{s}^{k_{s}})(A_{1}^{l_{1}}...A_{s}^{l_{s}})x=\mu_{1}^{k_{1}+l_{1}}...\mu_{s}^{k_{s}+l_{s}}x=(A_{1}^{l_{1}}...A_{s}^{l_{s}})(A_{1}^{k_{1}}...A_{s}^{k_{s}})x$$ 

hence $x\in\ker[A_{1}^{k_{1}}...A_{s}^{k_{s}},A_{1}^{l_{1}}...A_{s}^{l_{s}}]$ for any $k_{i},l_{j}\geq 0$. It follows that $x\in\CM(A_{1},...,A_{s})$ and $\CM(A_{1},...,A_{s})\neq 0$.

$\Leftarrow$ The subspace $\CM(A_{1},...,A_{s})$ is $A_{i}$-invariant by (1). Moreover, if $x\in\CM(A_{1},...,A_{s})$, then $x\in\ker [A_{i},A_{j}]$ for any $i,j=1,...,s$ and thus $A_{i}$ commute on $\CM(A_{1},...,A_{s})$. Since $\CM(A_{1},...,A_{s})\neq 0$, it follows by Corollary 2.3 that $A_{i}$ have a common eigenvector $v\in\CM(A_{1},...,A_{s})$. \epv

Now our aim is to show that the condition $$\CM(A_{1},...,A_{s})\neq 0$$ is computable. To prove this fact we need the following two lemmas.  

Assume that $A\in\MM_{n}(\CCC)$ and $f_{A}\in\CCC[x]$ is the characteristic polynomial of $A$. Recall that $f_{A}(A)=0$ by the renowned theorem of Cayley and Hamilton.

\begin{lem} Assume that $A\in\MM_{n}(\CCC)$ and $t\geq n$. Then the matrix $A^{t}$ is a linear combination of matrices $1_{n}=A^{0},A^{1},...,A^{n-1}$.
\end{lem}

{\bf Proof.} Assume that $f_{A}(x)=x^{n}+a_{n-1}x^{n-1}+...+a_{1}x+a_{0}$. Since $f_{A}(A)=0$, we have $A^{n}=-a_{n-1}A^{n-1}-...-a_{1}A-a_{0}1_{n}$. 

Assume that $t\geq n$ and all matrices $A^{n}, A^{n+1},..,A^{t}$ are linear combinations of $1_{n},A,...,A^{n-1}$. Then $A^{t+1}=A^{n}A^{t-n+1}=-a_{n-1}A^{t}-...-a_{1}A^{t-n+2}-a_{0}A^{t-n+1}$ and hence $A^{t+1}$ is a linear combination of $1_{n},A,...,A^{n-1}$. Thus the lemma follows by induction. \epv

We denote by $\langle\cdot | \cdot\rangle$ and $||\cdot||$ the standard scalar product in $\CCC^{n}$ and the standard norm in $\CCC^{n}$, respectively.

\begin{lem} Assume that $A,B\in\MM_{n}(\CCC)$. Then $$\ker (A^{*}A+B^{*}B)=\ker A \cap \ker B.$$
\end{lem}

{\bf Proof.} Assume that $X=A^{*}A$, $Y=B^{*}B$. Clearly $\ker A \cap \ker B\subseteq\ker (X+Y)$ so it is sufficient to show that $\ker (X+Y)\subseteq\ker A \cap \ker B$.

Assume that $\alpha\in\ker (X+Y)$. Then $X\alpha+Y\alpha=0$ and hence $(\ov{\alpha})^{tr}X\alpha+(\ov{\alpha})^{tr}Y\alpha=0$. Observe that for any $\beta\in\CCC^{n}$ we have $$(\ov{\beta})^{tr}X\beta=\langle\beta | X\beta\rangle=\langle\beta | A^{*}A\beta\rangle=\langle A\beta | A\beta\rangle=||A\beta||^{2}$$ and similarly $(\ov{\beta})^{tr}Y\beta=||B\beta||^{2}$. This implies that $X,Y$ are positive semi-definite and hence $(\ov{\alpha})^{tr}X\alpha+(\ov{\alpha})^{tr}Y\alpha=0$ yields that $$||A\alpha||^{2}=(\ov{\alpha})^{tr}X\alpha=0=(\ov{\alpha})^{tr}Y\alpha=||B\alpha||^{2}.$$

In consequence, $A\alpha=B\alpha=0$ and $\alpha\in \ker A \cap \ker B$. This shows that $\ker (X+Y)\subseteq\ker A \cap \ker B$.\epv

\begin{cor} Assume that $A_{1},...,A_{s}\in\MM_{n}(\CCC)$. Then the following equalities holds.

{\rm (1)} $$\CM(A_{1},...,A_{s})=\bigcap_{\begin{smallmatrix} k_{i},l_{j}\geq 0 \\ k_{1}+k_{2}+...+k_{s}\neq 0\\l_{1}+l_{2}+...+l_{s}\neq 0\end{smallmatrix}}^{n-1}\ker [A_{1}^{k_{1}}...A_{s}^{k_{s}},A_{1}^{l_{1}}...A_{s}^{l_{s}}].$$

{\rm (2)} $\CM(A_{1},...,A_{s})=\ker K$ where 

%$$L=\left[\begin{matrix}[A_{1}^{1}A_{2}^{0}...A_{s}^{0},A_{1}^{1}A_{2}^{0}...A_{s}^{0}]\\ [A_{1}^{1}A_{2}^{1}A_{3}^{0}...A_{s}^{0},A_{1}^{1}A_{2}^{0}...A_{s}^{0}] \\ \vdots \\ [A_{1}^{1}A_{2}^{1}...A_{s}^{1},A_{1}^{1}A_{2}^{1}...A_{s}^{1}] \\ [A_{1}^{2}A_{2}^{0}...A_{s}^{0},A_{1}^{1}A_{2}^{0}...A_{s}^{0}]\\\vdots \\ [A_{1}^{n-1}A_{2}^{n-1}...A_{s}^{n-1},A_{1}^{n-1}A_{2}^{n-1}...A_{s}^{n-1}]\end{matrix}\right],$$ 

$$K=\sum_{\begin{smallmatrix} k_{i},l_{j}\geq 0 \\ k_{1}+k_{2}+...+k_{s}\neq 0\\l_{1}+l_{2}+...+l_{s}\neq 0\end{smallmatrix}}^{n-1}[A_{1}^{k_{1}}...A_{s}^{k_{s}},A_{1}^{l_{1}}...A_{s}^{l_{s}}]^{*}[A_{1}^{k_{1}}...A_{s}^{k_{s}},A_{1}^{l_{1}}...A_{s}^{l_{s}}].$$

\end{cor}

{\bf Proof.} (1) Assume that $t_{1},...,t_{s}\geq 0$. It follows by Lemma 2.5 that $A_{1}^{t_{1}}...A_{s}^{t_{s}}$ is a linear combination of $A_{1}^{k_{1}}...A_{s}^{k_{s}}$ for some $0\leq k_{i}\leq n-1$. Moreover, observe that $$[A+B,C+D]=[A,C]+[A,D]+[B,C]+[B,D]$$ and hence $$\ker [A,C]\cap\ker [A,D]\cap\ker [B,C]\cap\ker [B,D]\subseteq\ker[A+B,C+D]$$ for any $A,B,C,D\in\MM_{n}(\CCC)$. This implies what is required.

(2) The equality follows by (1) and Lemma 2.6. \epv

\textit{Proof of Theorem 1.2.} (1) is a consequence of Theorem 2.4 (2) and Corollary 2.7 (1) whereas (2) of Corollary 2.7 (2).\epv

Note that Theorem 1.2 provide us a computable condition verifing the existence of a common eigenvector of $s\geq 2$ complex matrices $A_{1},...,A_{s}\in\MM_{n}(\CCC)$. The condition is a generalization of \cite[Theorem 3.1]{Sh} since $\CM(A,B)=\CN(A,B)$ for all $A,B\in\MM_{n}(\CCC)$, see Theorem 1.1. 

Indeed, it is easy to see from the proof of \cite[Theorem 3.1]{Sh} that if $v\in\CN(A,B)$, then $$(A^{p_{1}}B^{p_{2}})(A^{q_{1}}B^{q_{2}})v=(A^{p_{1}+q_{1}}B^{p_{2}+q_{2}})v=(A^{q_{1}}B^{q_{2}})(A^{p_{1}}B^{p_{2}})v$$ for all $p_{1},p_{2},q_{1},q_{2}\geq 0$ and hence $v\in\ker[A^{p_{1}}B^{p_{2}},A^{q_{1}}B^{q_{2}}]$. This implies that $\CN(A,B)\subseteq\CM(A,B)$ and since clearly $\CM(A,B)\subseteq\CN(A,B)$, we get $\CN(A,B)=\CM(A,B)$.

\section{Common invariant subspaces of higher dimensions}

In Section 2 we proved a computable criterion that allows one to check whether there exists a common eigenvector of $s\geq 2$ complex matrices $A_{1},...,A_{s}\in\MM_{n}(\CCC)$, or equivalently, a common invariant subspace of dimension $1$. 

In this section we apply the above result to show a computable criterion for existence of a common invariant subspace of dimension $d\geq 2$ of matrices $A_{1},...,A_{s}\in\MM_{n}(\CCC)$, if any $A_{i}$ has pairwise different eigenvalues. Our approach bases on methods from \cite{GI} and \cite{Ts}. 

Moreover, we show that the set of all $n\times n$ complex matrices having at least one multiple eigenvalue is Lebesgue-measurable, and of measure zero. This can can be interpret in the following way: if all matrices $A_{1},...,A_{s}\in\MM_{n}(\CCC)$ are random then it should be expected that each one of them has pairwise different eigenvalues. In consequence, the assumption seems not to be so strong in practical applications. 

First we recall the definition of \textit{$k$-th compound} of a matrix, see \cite{La}. 

Assume that $n\geq 1$ and $\alpha,\beta\subseteq\langle n\rangle$, where $\langle s\rangle:=\{1,...,s\}$ for any $s\in\NN$. If $A\in\MM_{n}(\CCC)$, we denote by $A[\alpha|\beta]$ a submatrix of $A$ composed from rows of $A$ indexed by $\alpha$ and columns of $A$ indexed by $\beta$. 

Assume that $k\leq n$ and $$Q_{k,n}:=\{(i_{1},...,i_{k})\in\NN^{k}|1\leq i_{1}<i_{2}<...<i_{k}\leq n\}$$ is the set of all $k$-tuples of elements from $\langle n\rangle$ ordered lexicographically.

The \textit{$k$-th compound} of $A\in\MM_{n}(\CCC)$ is defined to be the matrix $$C_{k}(A):=[\det A[\alpha|\beta]]_{\alpha,\beta\in Q_{k,n}}\in\MM_{{n\choose k}}(\CCC).$$

The following theorem is implicitly contained in \cite{Ts}. It generalizes \cite[Theorem 2.2]{GI} and \cite[Theorem 3.1]{GI} for the case of $s$ complex matrices where $s\geq 2$.

\begin{thm} Assume that $A_{1},...,A_{s}\in\MM_{n}(\CCC)$ are nonsingular, $k\in\NN$, $2\leq k<n$ and $C_{k}(A_{1}),...,C_{k}(A_{s})$ have pairwise different eigenvalues. Then $A_{i}$ have a common invariant subspace of dimension $k$ if and only if $C_{k}(A_{i})$ have a common eigenvector.
\end{thm}

{\bf Proof.} Theorem follows easily from the proof of \cite[Theorem 2.2]{Ts}. \epv

We can exchange assumptions of Theorem 3.1 for the assumption that $A_{1},...,A_{s}$ have pairwise different eigenvalues. It can be easily concluded from the proposition below which is essentially the same as \cite[Lemma 2.4]{GI}.

\begin{prop} Assume that $A\in\MM_{n}(\CCC)$ has pairwise different eigenvalues, $k\in\NN$ and $2\leq k<n$. Then there is $t\in\{0,...,p\}$ where $$p=\frac{{n\choose k}({n\choose k}-1)}{2}(k-1)$$ such that $A-t1_{n}$ is nonsingular and $C_{k}(A-t1_{n})$ has pairwise different eigenvalues.
\end{prop}

{\bf Proof.} It follows by \cite[Lemma 2.4]{GI} that there is $t\in\{0,...,p\}$ where $$p=\frac{{n\choose k}({n\choose k}-1)}{2}(k-1)$$ such that $C_{k}(A-t1_{n})$ has pairwise different eigenvalues. All eigenvalues of $C_{k}(A-t1_{n})$ are of the form $(t-\lambda_{i_{1}})(t-\lambda_{i_{2}})...(t-\lambda_{i_{k}})$ where $1\leq i_{1}<i_{2}<...<i_{k}\leq n$ and $\lambda_{1},\lambda_{2},...,\lambda_{n}$ are eigenvalues of $A$, see \cite[Theorem 2.1 (7)]{GI}. This implies that $t$ is not an eigenvalue of $A$ since otherwise it is clear that $0$ is a multiple eigenvalue of $C_{k}(A-t1_{n})$. Hence $\det(A-t1_{n})\neq 0$ and $A-t1_{n}$ is nonsingular. \epv

\begin{cor} Assume that $A_{1},...,A_{s}\in\MM_{n}(\CCC)$ have pairwise different eigenvalues and $k\in\NN$, $2\leq k<n$.

{\rm (1)} There are $t_{1},...,t_{s}\in\{0,...,p\}$ where $$p=\frac{{n\choose k}({n\choose k}-1)}{2}(k-1)$$ such that $\wt{A_{i}}:=A_{i}-t_{i}1_{n}$ are nonsingular and $C_{k}(\wt{A_{i}})$ have pairwise different eigenvalues. 

{\rm (2)} The matrices $A_{1},...,A_{s}\in\MM_{n}(\CCC)$ have a common invariant subspace of dimension $k$ if and only if $C_{k}(\wt{A_{i}})$ have a common eigenvector.

{\rm (3)} The matrices $A_{1},...,A_{s}\in\MM_{n}(\CCC)$ have a common invariant subspace of dimension $k$ if and only if $\ker K\neq 0$ where

%$$L=\left[\begin{matrix}[X_{1}^{1}X_{2}^{0}...X_{s}^{0},X_{1}^{1}X_{2}^{0}...X_{s}^{0}]\\ [X_{1}^{1}X_{2}^{1}X_{3}^{0}...X_{s}^{0},X_{1}^{1}X_{2}^{0}...X_{s}^{0}] \\ \vdots \\ [X_{1}^{1}X_{2}^{1}...X_{s}^{1},X_{1}^{1}X_{2}^{1}...X_{s}^{1}] \\ [X_{1}^{2}X_{2}^{0}...X_{s}^{0},X_{1}^{1}X_{2}^{0}...X_{s}^{0}]\\\vdots \\ [X_{1}^{n-1}X_{2}^{n-1}...X_{s}^{n-1},X_{1}^{n-1}X_{2}^{n-1}...X_{s}^{n-1}]\end{matrix}\right],$$

$$K=\sum_{\begin{smallmatrix} k_{i},l_{j}\geq 0 \\ k_{1}+k_{2}+...+k_{s}\neq 0\\l_{1}+l_{2}+...+l_{s}\neq 0\end{smallmatrix}}^{n-1}[X_{1}^{k_{1}}...X_{s}^{k_{s}},X_{1}^{l_{1}}...X_{s}^{l_{s}}]^{*}[X_{1}^{k_{1}}...X_{s}^{k_{s}},X_{1}^{l_{1}}...X_{s}^{l_{s}}]$$ and $X_{i}:=C_{k}(\wt{A_{i}})$ for $i=1,...,s$.

\end{cor}

{\bf Proof.} (1) follows by Proposition 3.2. (2) follows by (1), Theorem 3.1 and a simple fact that $A_{i}$ have a common invariant subspace $W$ if and only if $A_{i}-q_{i}1_{n}$ have a common invariant subspace $W$, for any $W\subseteq\CCC^{n}$ and $q_{i}\in\CCC$. (3) follows by (2) and Theorem 1.2.\epv

Now we propose a finite and deterministic algorithm verifing the existence of a common invariant subspace of $A_{1},...,A_{s}\in\MM_{n}(\CCC)$ of a fixed dimension $d$ provided every $A_{i}$ has pairwise different eigenvalues.

Recall first that the \textit{discriminant} $\disc(f)$ of a polynomial $f\in\CCC[x]$ is, by the definition, the resultant of $f$ and $f'$ where $f'$ denotes the formal derivative of $f$, see \cite[Chapter IV]{Lang}, Section 8. 

It is commonly known that $\disc(f)=0$ if and only if $f$ has a multiple root. Obviously, the condition $\disc(f)=0$ is computable.\\

{\bf Algorithm.} Input: $A_{1},...,A_{s}\in\MM_{n}(\CCC)$ having pairwise different eigenvalues and $d\in\{1,...,n-1\}$. Output: 'yes' if there is a common invariant subspace of dimension $d$ and 'no' otherwise.
\begin{enumerate}
	\item If $d=1$, compute $\CM(A_{1},...,A_{s})$ and check whether $\CM(A_{1},...,A_{s})\neq 0$ using Theorem 1.2. If $\CM(A_{1},...,A_{s})\neq 0$, print 'yes', otherwise print 'no'.
	\item If $d\geq 2$, compute $p=\frac{{n\choose d}({n\choose d}-1)}{2}(d-1)$ and $t_{i}\in\{0,...,p\}$ such that $\det(\wt{A_{i}})\neq 0$ and $\disc(f_{i})\neq 0$ where $\wt{A_{i}}:=A-t_{i}1_{n}$ and $f_{i}$ denotes the characteristic polynomial of $C_{d}(\wt{A_{i}})$. Go to the step 3.
	\item Compute $\ker K$ where the matrix $K$ is the matrix from Corollary 3.3 (3). If $\ker K\neq 0$, print 'yes', otherwise print 'no'.\epv
\end{enumerate}

The correctness of the above algorithm follows from Theorem 1.2 and Corollary 3.3. 

Our next aim is to show that the set of all complex $n\times n$ matrices having at least one multiple eigenvalue is Lebesgue-measurable, and of measure zero. We set the following notation.

We identify the linear space $\MM_{n}(\CCC)$ of all complex $n\times n$ matrices with $\CCC^{n^{2}}$ which is a measurable space with the Lebesgue measure induced from $\RR^{2n^{2}}$, see \cite[Chapter I]{GR} for details. Hence the space $\MM_{n}(\CCC)$ becomes a measurable space with the same Lebesgue measure as the one on $\RR^{2n^{2}}$. We denote this measure by $\CL_{n}$.

Assume that $A=[a_{ij}]_{i,j=1,...,n}\in\MM_{n}(\CCC)$ and $f_{A}(t)=\det(A-t1_{n})\in\CCC[t]$ is the characteristic polynomial of $A$. Although the discriminant $\disc(f_{A})$ of $f_{A}$ is a complex number, it is clear that it can be viewed as a polynomial of $n^{2}$ variables $a_{ij}$ with integer coefficients.

Hence we define $\disc_{n}\in\ZZ[(x_{ij})_{i,j=1,...,n}]$ in such a way that $$\disc_{n}((a_{ij})_{i,j=1,...,n})=\disc(f_{A})$$ for any $A=[a_{ij}]_{i,j=1,...,n}\in\MM_{n}(\CCC)$.

Observe that the set $$\CZ_{n}:=\{A=[a_{ij}]_{i,j=1,...,n}\in\CCC^{n^{2}}|\disc_{n}((a_{ij})_{i,j=1,...,n})=0\}$$ represents the set of all complex $n\times n$ matrices having at least one multiple eigenvalue.

\begin{prop} The set $\CZ_{n}$ is Lebesgue-measurable and $\CL_{n}(\CZ_{n})=0$.
\end{prop}

{\bf Proof.} Clearly, $\disc_{n}$ viewed as a complex function is measurable and holomorphic in $\CCC^{n^{2}}$. Moreover, $\disc_{n}$ is not identically zero since there are matrices with pairwise different eigenvalues. Hence \cite[Chapter I, Corollary 10]{GR} yields $\CL_{n}(\CZ_{n})=0$.\epv

\section{An application to completely positive superoperators}

This section is devoted to show an application of our results in the theory of completely positive superoperators defined on a finite-dimensional complex Hilbert space, see \cite{HZ}, \cite{Ja}. 

Assume that $\CH$ is a complex Hilbert space and $B(\CH)$ is the set of all linear continouous operators on $\CH$. It is commonly known that $B(\CH)$ is equipped with the structure of a Banach space. As it was mentioned in Section 1, a linear and continouous map $\Phi:B(\CH)\ra B(\CH)$ on $B(\CH)$ is called a \textit{superoperator}.

Assume that $\CH$ is a complex finite-dimensional Hilbert space and $\dim(\CH)=n$. A superoperator $\Phi:B(\CH)\ra B(\CH)$ on $B(\CH)$ of the form $$\Phi(X)=\sum_{i=1}^{s}K_{i}XK_{i}^{*},$$ where $K_{i}\in B(\CH)$ is \textit{completely positive}. Such operators play a prominent role in quantum information theory, see \cite{HZ}, \cite{Ja} for details.

The operators $K_{1},...,K_{s}\in B(\CH)$ in the above formula are called the \textit{Kraus coefficients} of $\Phi$. Clearly, they can be viewed as arbitrary $n\times n$ complex matrices and also can be treated as objects that define $\Phi$.

Important subclass of the class of all completely positive superoperators is formed by \textit{irreducible} completely positive superoperators, see \cite{Fa} for the definitions. The following theorem connects irreducible completely positive superoperators with the subject matter of the paper. 

\begin{thm} \textnormal{(see \cite{Fa})} Assume that $\CH$ is a complex finite-dimensional Hilbert space and $\Phi:B(\CH)\ra B(\CH)$ is a completely positive superoperator on $B(\CH)$ such that $\Phi(X)=\sum_{i=1}^{s}K_{i}XK_{i}^{*}$. Then $\Phi$ is irreducible if and only if the matrices $K_{i}$ do not have a nontrivial common invariant subspace in $\CH$.\epv
\end{thm}

In view of Theorem 4.1 we can apply the algorithm presented in Section 3 to check whether given completely positive superoperator $\Phi$ is irreducible, if its Kraus coefficients do not have a multiple eigenvalue. Moreover, Proposition 3.4 yields that such a situation should be quite common in practical applications.

\noindent 
Andrzej Jamio{\l}kowski\\
Faculty of Physics, Astronomy and Informatics\\ 
Nicolaus Copernicus University\\
Grudzi\k{a}dzka 5\\
87-100 Toru\'n, Poland\\
jam@fizyka.umk.pl\\\\
Grzegorz Pastuszak\\
Center for Theoretical Physics of the Polish Academy of Sciences\\
Al. Lotnik{\'o}w 32/46\\
02-668 Warszawa, Poland\\
past@mat.umk.pl

\end{document}